\documentclass[10pt,reqno]{amsart}
\usepackage{mathrsfs}
\usepackage{amssymb}
\usepackage{amsmath}
\usepackage{amsfonts}
\usepackage{amscd}
\usepackage{multirow}
\usepackage{graphicx}
\usepackage{tabularx}
\usepackage{longtable}
\usepackage[colorlinks, linkcolor=blue, anchorcolor=blue, citecolor=blue]{hyperref}
\begin{document}
\setlength{\oddsidemargin}{0cm} \setlength{\evensidemargin}{0cm}
\baselineskip=22pt

\theoremstyle{plain} \makeatletter
\newtheorem{theorem}{Theorem}[section]
\newtheorem{proposition}[theorem]{Proposition}
\newtheorem{lemma}[theorem]{Lemma}
\newtheorem{coro}[theorem]{Corollary}

\theoremstyle{definition}
\newtheorem{defi}[theorem]{Definition}
\newtheorem{notation}[theorem]{Notation}
\newtheorem{exam}[theorem]{Example}
\newtheorem{prop}[theorem]{Proposition}
\newtheorem{conj}[theorem]{Conjecture}
\newtheorem{prob}[theorem]{Problem}
\newtheorem{remark}[theorem]{Remark}
\newtheorem{claim}{Claim}

\newcommand{\SO}{{\mathrm S}{\mathrm O}}
\newcommand{\SU}{{\mathrm S}{\mathrm U}}
\newcommand{\U}{{\mathrm U}}
\newcommand{\C}{{\mathbb C}}
\newcommand{\Sp}{{\mathrm S}{\mathrm p}}
\newcommand{\so}{{\mathfrak s}{\mathfrak o}}
\newcommand{\Ad}{{\mathrm A}{\mathrm d}}
\newcommand{\ad}{{\mathrm a}{\mathrm d}}
\newcommand{\m}{{\mathfrak m}}
\newcommand{\n}{{\mathfrak n}}
\newcommand{\g}{{\mathfrak g}}
\newcommand{\p}{{\mathfrak p}}
\newcommand{\q}{{\mathfrak q}}
\newcommand{\fk}{{\mathfrak k}}
\newcommand{\fh}{{\mathfrak h}}
\newcommand{\E}{{\mathrm E}}
\newcommand{\F}{{\mathrm F}}
\newcommand{\ffb}{{\mathfrak b}}
\newcommand{\fl}{{\mathfrak l}}
\newcommand{\ft}{{\mathfrak t}}
\newcommand{\G}{{\mathrm G}}
\newcommand{\K}{{\bar{K}}}
\newcommand{\R}{\mathrm{Ric}}
\newcommand{\expo}{\mathrm{exp}}
\newcommand{\yd}{\approx}
\newcommand{\zz}{\mathbb{Z}_2}
\newcommand{\zzzz}{\mathbb{Z}_4}

\numberwithin{equation}{section}

\newcolumntype{Y}{>{\centering\arraybackslash}X}

\title[Geodesic orbit metrics on certain homogeneous spaces]{Geodesic orbit metrics on homogeneous spaces constructed by strongly isotropy irreducible spaces}
\author{Huibin Chen, Zhiqi Chen and Fuhai Zhu$^\dagger$}
\address{School of Mathematical Sciences, Nanjing Normal University, Nanjing, 210023, P.R. China}
\email{chenhuibin@njnu.edu.cn}
\address{School of Mathematical Sciences and LPMC, Nankai University, Tianjin 300071, P.R. China}
\email{chenzhiqi@nankai.edu.cn}
\address{Department of Mathematics, Nanjing University, Nanjing 210023, P.R. China}
\email{zhufuhai@nju.edu.cn}
\date{}
\thanks{$^\dagger$Corresponding Author}
\date{}

\maketitle
\begin{abstract}
In this paper, we focus on homogeneous spaces which are constructed from two strongly isotropy irreducible spaces, and prove that any geodesic orbit metric on these spaces is naturally reductive.
\end{abstract}

\section{Introduction}
It is well-known that any geodesic in a Riemannian symmetric space $G/K$ is an orbit of some one-parameter subgroup of $G$. As a natural generalization of Riemannian symmetric spaces, O. Kowalski and L. Vanhecke introduced the notion of a geodesic orbit space in \cite{KoVa}. A homogeneous Riemannian manifold $(G/K, g)$ is called a geodesic orbit space if any geodesic of $G/K$ is an orbit of a one-parameter subgroup of $G$ and $g$ is called a geodesic orbit metric. For more information on geodesic orbit spaces, one can refer to the survey \cite{Ar}. For left-invariant geodesic orbit metrics on Lie groups, one can refer to \cite{ChChDe, ChChWo} and \cite{Ni}. For the generalizations to pseudo-Riemannian manifolds and Finsler manifolds, see \cite{Ba, ChWo, Du, La, YaDe} and the reference therein.

Naturally reductive spaces form an important class of geodesic orbit spaces. In \cite{KoVa}, O. Kowalski and L. Vanhecke proved that every geodesic orbit space of dimension $\leq 5$ is naturally reductive. In \cite{Ka}, A. Kaplan constructed the first example of non-naturally reductive geodesic orbit space.

Weakly symmetric spaces, introduced by A. Selberg in \cite{Se}, form another important class of geodesic orbit spaces, as is proved by J. Berndt, O. Kowalski and L. Vanhecke in \cite{BeKoVa}. Riemannian weakly symmetric spaces with reductive isometry groups were classified by O.S. Yakimova in \cite{Ya}.

In \cite{Go}, C.S. Gordon reduced the classification problem of geodesic orbit manifolds to three special cases, but it seems still impossible to classify all geodesic orbit manifolds. Thus many scholars turned their attention to construct examples under some special settings, for example, the results in \cite{AlAr, AlNi, ArWa, ArWaZh, ArWaZh2, ChNi, DuKoNi, Go, KoVa, NiNi, Ta} and so on. 

For a given homogeneous space $G/H$, it is a fundamental problem to determine all Riemannian geodesic orbit metrics. In this paper, we study a class of irreducible compact homogeneous spaces $G/H$ of semisimple type which are not weakly symmetric spaces, and determine all $G$-invariant geodesic orbit metrics up to homothety on each of these homogeneous spaces.

One of the main tools used here is the principal isotropy subgroups of irreducible representations of simple Lie groups, for more details see \cite{El, HiHi, Kr}. The relation between the structure of Riemannian geodesic orbit spaces and the principal isotropy subgroups of representations of Lie groups was first noticed and studied by Yu.G. Nikonorov in \cite{Ni2}. Another tool is so-called geodesic graph, which was originally introduced in \cite{KoVa} to prove when a geodesic orbit metric is naturally reductive. The main theorem of this paper is given as follows.
\begin{theorem}\label{mainthm}
Let $G_i/K$ $(i=1,2)$ be strongly isotropy irreducible spaces with $K$ simple. Then any geodesic orbit metric on $(G_1\times G_2)/\Delta K$ is naturally reductive.
\end{theorem}

The paper is organized as follows. In section 2, we will recall some basic facts on geodesic orbit manifolds and some equivalent conditions for a Riemannian homogeneous space to be a geodesic orbit space. In section 3, we will introduce homogeneous spaces $(G_1\times G_2)/\Delta K$ arising from two strongly isotropy irreducible spaces $G_1/K$ and $G_2/K$. The classification when $G_1\not=G_2$ is presented in Appendix I. In section 4, we will characterize principal isotropy subgroups of irreducible representations appeared in Appendix II. In section 5, we will determine all geodesic orbit metrics on $(G_1\times G_2)/\Delta K$ for $G_1\not=G_2$. In section 6, by a different method, we will classify geodesic orbit metrics on $(G_1\times G_1)/\Delta K$ where $G_1/K$ is a strongly isotropy irreducible space.

\section{Geodesic orbit manifolds}
Let $(G/K, g)$ be a Riemannian homogeneous space with $G$ compact, connected and semisimple. Denote by $\g$ and $\fk$ the corresponding Lie algebras of $G$ and $K$ respectively.
A geodesic $\gamma(t)$ on $G/K$ through the origin $o=eK$ is called homogeneous if it is an orbit of a one-parameter subgroup of $G$, i.e.
\begin{equation}
\gamma(t)=\expo(tX)\cdot o,\quad t\in\mathbb{R},
\end{equation}
where $X$ is a non-zero vector in $\g$ and we call $X$ a geodesic vector.

Let $B$ be the minus Killing from of $\g$. Consider the $B$-orthogonal $reductive$ $decomposition$ of $\g$:
\begin{equation}\label{rd}
\g=\fk\oplus\m.
\end{equation}
Then we have the following result (Geodesic Lemma).
\begin{lemma}[Proposition 2.1, \cite{KoVa}]
A nonzero vector $X\in\g$ is a geodesic vector if and only if for all $Y\in\m$
\begin{equation*}
([X, Y ]_\m , X_\m)=0,
\end{equation*}
where the subscript $\m$ denotes the projection into $\m$.
\end{lemma}

\begin{defi}
A homogeneous Riemannian manifold $(G/K, g)$ is called a geodesic orbit space (or g.o. space) if all geodesics on $G/K$ are homogeneous with respect to $G$. In this case, the metric $g$ is called a geodesic orbit metric (or g.o. metric). If $G$ is the full isometry group, then $(G/K, g)$ is called a geodesic orbit manifold (or g.o. manifold).
\end{defi}

Any $G$-invariant Riemannian metric on $G/K$ determines an $\Ad(K)$-invariant inner product $(\ ,\ )$ on $\m$, which in turn determines an $\Ad(K)$-equivariant positive definite symmetric operator $A$ on $\m$ satisfying $(x, y)=B(Ax,y)$ for all $x,y\in\m$, and vice versa. The operator $A$ is called a metric endomorphism on $\m$. A useful description of homogeneous geodesics is the following theorem.
\begin{theorem}[Proposition 1, \cite{AlAr}]
Let $(G/K, g)$ be a homogeneous Riemannian manifold and $A$ the corresponding metric endomorphism. Let $Z\in\fk$ and $X\in\m$. Then the following are equivalent:
\begin{enumerate}
\item $Z+X$ is a geodesic vector in $\g$.
\item $[Z+X, A(X)]\in\fk$.
\item $([Z, X], Y)=(X, [X, Y]_\m)$ for all $Y\in\m$.
\item $[Z+X, Y]_\m, X)=0$ for all $Y\in\m$.
\end{enumerate}
\end{theorem}
As a consequence, we obtain the following characterization of geodesic orbit spaces.

\begin{coro}[Corollary 2, \cite{AlAr}]\label{GOcon}
Let $(G/K, g)$ be a homogeneous Riemannian manifold. Then $(G/K, g)$ is a geodesic orbit space if and only if for every $X\in\m$ there exists some $Z\in\fk$ such that
\begin{equation} \label{shi123}
[Z+X, A(X)]\in\fk,
\end{equation}
where $A$ is the corresponding metric endomorphism on $\m$.
\end{coro}

In fact, this corollary could be strengthened to the following lemma:
\begin{lemma}[Proposition 2, \cite{So}]\label{lemma=0}
A Riemannian homogenous space $(G/K, g)$ is geodesic orbit if and only if for any $X\in\m$ there exists $Z\in\fk$ such that
\begin{equation}
[Z+X, AX]=0.
\end{equation} 
\begin{proof}
By Corollary \ref{GOcon}, we only have to prove $[X, AX]\in\m$ because $[Z, AX]\in\m$ always holds for any $Z\in\fk$. Since the metric endomorphism $A$ is $\ad(\fk)$-equivariant and the Killing form $B$ is $\Ad(H)$-invariant and $A$-symmetric, we have for any $\tilde{Z}\in\fk$,
\begin{equation}
B([X, AX],\tilde{Z})=-B(X, [\tilde{Z}, AX])=-B(X, A[\tilde{Z}, X])=-B(AX, [\tilde{Z}, X])=-B([X, AX], \tilde{Z}),
\end{equation}
which means $B([X, AX], \tilde{Z})$ holds for any $\tilde{Z}\in\fk$. Therefore, $[X, AX]\in\m$. This completes the proof of this theorem.
\end{proof}
\end{lemma}

A homogeneous Riemannian manifold $(G/K, g)$ is called naturally reductive if there is an $\Ad(K)$-invariant decomposition $\g=\fk+\m$ such that for all $X,Y,Z\in\m$,
\begin{equation}\label{nri}
([X, Y]_\m, Z)+(Y, [X, Z]_\m)=0.
\end{equation}
An equivalent expression of the identity (\ref{nri}) is that for all $X,Y\in\m$,
\begin{equation}
([X, Y]_\m, X)=0.
\end{equation}
It is well-known that natural reductivity is equivalent to the geometrical property that for each vector $X\in\m$, the orbit $\gamma(t)=\expo tX\cdot o$ is a geodesic. As a result, naturally reductive Riemannian homogeneous spaces are geodesic orbit spaces, which can be also deduced from the algebraic equivalent conditions.

It is notable that the natural reductivity is adapted to the $\Ad(K)$-invariant decomposition $\g=\fk+\m$. This property leads to an obstacle to investigate whether a geodesic orbit metric is naturally reductive, in fact, the $\Ad(K)$-invariant decomposition is not unique when the isometry group $G$ is non-simple. The following lemma provides an effective method to determine when a geodesic orbit metric is naturally reductive.
\begin{lemma}[Proposition 2.10, \cite{KoVa}]\label{geogra}
	Let $G/K$ be a Riemannian homogeneous space with a given $\Ad(K)$-invariant decomposition $\g=\fk+\m$. Then $G/K$ is naturally reductive if and only if there is an $\Ad(K)$-equivariant linear map $\xi:\m\rightarrow\fk$ such that $Z=\xi(X)$ satisfying relation (\ref{shi123}).
\end{lemma}

\section{A class of compact semisimple homogeneous spaces}

J.A. Wolf \cite{Wo} gave a complete classification of coset spaces $G/K$ satisfying the following conditions(also see O.V. Manturov \cite{Ma, Ma2, Ma3} and M. Kr\"{a}mer \cite{Kr2}):
\begin{enumerate}
\item $G$ is a connected Lie group and $K$ is a closed subgroup,
\item $G/K$ is a reductive coset space on which $G$ acts effectively,
\item the linear isotropy action (on the tangent space of $G/K$) of the identity component $K_0$ of $K$ is a representation which is irreducible over the real number field,
\item $K$ is compact,
\item $(G, K)$ is not a symmetric pair, i.e. $K_0$ is not the identity component of the fixed point set of an involution of $G$.
\end{enumerate}

These spaces were renamed as strongly isotropy irreducible spaces by M. Wang and W. Ziller in \cite{WaZi}. Let $G_1/K$ and $G_2/K$ be two strongly isotropy irreducible spaces of compact type with $K$ connected and simple. By the classification of the isotropy representations of non-symmetric strongly isotropy irreducible spaces given in the tables on pages 107-110 (Theorem 11.1 \cite{Wo}) and the classification of the isotropy representations of irreducible symmetric spaces given in tables (8.11.2) and (8.11.5) of \cite{Wo2}, it is not hard to give the classification of compact homogeneous space $G/H=(G_1\times G_2)/\Delta K$. The classification for $G_1\not=G_2$ is given in Appendix I, also we will give some information about the isotropy representations which is useful for classifying geodesic orbit metrics.
   
\begin{remark}
Here we get the classification directly from the group level, which is different from the classification given in \cite{YaChDe} by Z. Yan, H. Chen and S. Deng. In fact, for any two different strongly isotropy irreducible spaces $G_i/K_i$ $(i=1,2)$ with the Lie algebras $\fk_i$ of $K_i$ $(i=1,2)$ coincide, they got the classification of $(\mathfrak g_1\oplus\mathfrak g_2,\Delta\mathfrak k)$, and then obtained the classification of $(G_1\times G_2)/\widetilde{\Delta K}$ where $\widetilde{\Delta K}$ is the connected Lie subgroup of $G_1\times G_2$ with the Lie algebra $\Delta\mathfrak k$. Here, $K_1$ may not be isomorphic to $K_2$. While in our classification, we require $K_1\cong K_2$.

%In fact, they first gave the classification of the pair $(\mathfrak g_i,\mathfrak k)$ for the strongly isotropy irreducible spaces $G_i/K$ $(i=1,2)$, then got the classification of $(\mathfrak g_1\oplus\mathfrak g_2,\Delta\mathfrak k)$ in the Lie algebraic level, finally obtained the classification of $G_1\times G_2/\widetilde{\Delta K}$ where $\widetilde{\Delta K}$ is the connected Lie subgroup of $G_1\times G_2$ with the Lie algebra $\Delta\mathfrak k$. Here, $\widetilde{\Delta K}$ may not be the diagonal embedding of $K$ into $G_1\times G_2$. In our classification, we require $\Delta K$ to be the diagonal embedding of $K$ into $G_1\times G_2$.

%The difference is that we eliminate the following case: the isotropy irreducible spaces $G_1/K_1$ and $G_2/K_2$, $K_1\not=K_2$, but the Lie algebras of $K_1$ and $K_2$ are the same.
\end{remark}

\section{Principal isotropy subgroups of representations of compact simple Lie groups}
This section is to study principal isotropy subgroups of representations appeared in Appendix I. First we recall the study on principal isotropy subgroups of representations of connected compact simple Lie groups, which plays a fundamental role in next section.

Let $K$ be a compact linear Lie group acting linearly on a finite-dimensional vector space $V$. Then almost all points of $V$ are situated on the orbits of $K$, which are pairwise isomorphic as $K$-manifolds. These orbits are called $orbits$ $in$ $general$ $position$. The isotropy subgroups of all points in all orbits of general position are conjugated in $K$. Their conjugation class is called $a$ $principal$ $isotropy$ $subgroup$ for the linear Lie group $K$ and the corresponding Lie algebra is called $a$ $principal$ $isotropy$ $subalgebra$ or $a$ $stationary$ $subalgebra$ of points in general position (see \cite{ChNi, El, HiHi} and the reference therein).

\begin{defi}[\cite{HiHi}]
Let $\rho$ be a real irreducible representation of a simple Lie group $K$. We say that a principal isotropy subgroup $(H_\rho)$ ($H_\rho$ is a representative in the conjugacy class $(H_\rho)$) is trivial if
$H_\rho=\mathrm{Ker}(\rho)\subset Z(K),$
where $Z(K)$ is the center of $K$.
\end{defi}

The classification of representations of connected compact simple Lie groups with non-trivial connected principal isotropy subgroups has been given in \cite{El, HiHi} and \cite{Kr}. For readers' convenience, we will list the classification results in Theorem~\ref{HiHi} in Appendix II.

With the contents in Theorem~\ref{HiHi}(see Appendix II), we can determine the principal isotropy subgroups $(H_{\rho_i})$ of $\rho_i$ $(i=1,2)$ for each homogeneous space $G/H=(G_1\times G_2)/\Delta K$ with $G_1\neq G_2$ listed in Appendix I. For example, in case D.1, $\Delta K=\Delta\mathrm{G}_2$, $G_1\times G_2=\mathrm{Spin}(7)\times\mathrm{Spin}(14)$, $\rho_1=\varphi_1$ and $\rho_2=3\varphi_1$. Then we have $(H_{\rho_1})$ is a maximal torus and $(H_{\rho_2})$ is trivial from Appendix II. With a detailed case by case check, we can obtain the following theorem.

\begin{theorem}\label{pis}
Let the notation be as above. For every homogeneous space $G/H=(G_1\times G_2)/\Delta K$ for $G_1\not=G_2$ in Appendix I, one of the following conditions is satisfied:
\begin{enumerate}
 \item either the principal isotropy subalgebra $(\fh_{\rho_1})$ or $(\fh_{\rho_2})$ is $\{0\}$,
 \item either $(\fh_{\rho_1})$ or $(\fh_{\rho_2})$ is a torus.
\end{enumerate}
\end{theorem}

\section{Geodesic orbit metrics on $(G_1\times G_2)/\Delta K$ for $G_1\not=G_2$}
Let $G/H=(G_1\times G_2)/\Delta K$ be a homogenous space constructed as above, where $G_i/K(i=1,2)$ are strongly isotropy irreducible spaces with $K$ simple and connected. Denote the corresponding Lie algebras by $\mathfrak{g}_1$, $\mathfrak{g}_2$ and $\mathfrak{k}$, respectively. Let $B_1$, $B_2$ and $B_\mathfrak{k}$ be the corresponding minus Killing form of $\frak{g}_1$, $\frak{g}_2$ and $\frak{k}$. Then there are two constants $c_1, c_2>0$ such that $B_\frak{k}=c_1B_1|_{\fk}$ and $B_\frak{k}=c_2B_2|_{\fk}$. Let $B=B_1+B_2$ be the minus Killing form of $\frak{g}=\frak{g}_1\oplus\frak{g}_2$. Let $\g_i=\fk\oplus\p_i$ be the $\Ad(K)$-invariant decomposition of $\g_i$ with $B_i(\fk, \p_i)=0$, $i=1,2$. Then we have the corresponding $B$-orthogonal $reductive$ $decomposition$ of $\frak{g}$:
\begin{equation}\label{dc1}
\g=\fh\oplus\m=\fh\oplus\m_0\oplus\m_1\oplus\m_2,
\end{equation}
where $\fh=\{(Z, Z)|Z\in\fk\}$, $\m_0=\{(Z, -\frac{c_2}{c_1}Z)|Z\in\fk\}$, $\m_1=\{(P, 0)|P\in\p_1\}$ and $\m_2=\{(0, P)|P\in\p_2\}$.

Here $G_i/K$ $(i=1,2)$ are different  strongly isotropy irreducible spaces and $K$ is simple, $\m_0, \m_1$ and $\m_2$ are irreducible and mutually non-equivalent as $\ad(\fk)$-modules. Hence, any $G$-invariant metric $g$ on $G/H$ can be determined by an $\Ad(H)$-invariant inner product on tangent space $T_o(G/H)\cong\m=\m_0\oplus\m_1\oplus\m_2$:
\begin{equation}\label{mt1}
(\ ,\ )'=x_0B|_{\m_0}+x_1B|_{\m_1}+y_1B|_{\m_2},
\end{equation}
for some $x_0, x_1,y_1\in \mathbb{R}^+$. Clearly, for any $a\in\mathbb{R}^+$, $(\ ,\ )'$ is a geodesic orbit metric if and only if $a(\ ,\ )'$ is a geodesic orbit metric. Thus, it is enough to study $G$-invariant metrics of the form
\begin{equation}\label{metric}
(\ ,\ )=B|_{\m_0}+xB|_{\m_1}+yB|_{\m_2},
\end{equation}
for some $x,y\in \mathbb{R}^+$. 

%First, we need to determine when the metric (\ref{metric}) on $G/H$ is naturally reductive.
%\begin{lemma}\label{nr}
%The $G$-invariant metric $g$ on $G/H$ induced by the inner product $(\ref{metric})$ is naturally reductive if and only if $x=y=1$.
%\end{lemma}
%\begin{proof}
%If the $G$-invariant metric $g$ on $G/H$ is naturally reductive, then for any $X,Y\in\m$,
%\begin{equation}
%0=([X, Y]_\m, X)=B([X, Y]_\m, AX)=B([X, Y], AX)=B(Y, [AX, X]),
%\end{equation}
%where $A$ is the metric endomorphism. Since $Y$ is arbitrary, we have $[AX, X]\in\fh$ for any $X\in\m$.
%For any $X=X_0+X_1+X_2\in\m$, $X_0=(K_0, -\frac{c_2}{\alpha c_1}K_0)\in\m_0, X_1=(P_1, 0)\in\m_1, X_2=(0, P_2)\in\m_2$, we have
%\begin{align*}
%[AX, X]&=[(K_0+xP_1,-\frac{c_2}{\alpha c_1}K_0+yP_2), (K_0+P_1,-\frac{c_2}{\alpha c_1}K_0+P_2)]\\
%&=([K_0+xP_1, K_0+P_1], [-\frac{c_2}{\alpha c_1}K_0+yP_2,-\frac{c_2}{\alpha c_1}K_0+P_2])\\
%&=((1-x)[K_0, P_1], \frac{c_2}{\alpha c_1}(y-1)[K_0, P_2])\in\fh.
%\end{align*}
%Since $K_0,P_1$ and $P_2$ are arbitrary, we have $x=y=1$. The converse is obvious.
%\end{proof}

First, we will deduce some equivalent conditions for Riemannian homogeneous spaces $(G/H, g)$ defined above to be geodesic orbit spaces.
\begin{lemma}\label{ec}
The $G$-invariant metric $g$ on $G/H$ induced by the inner product $(\ref{metric})$ is geodesic orbit if and only if for any $K_0\in\fk, P_1\in\p_1$ and $P_2\in\p_2$, there exists some $Z\in\fk$ such that
\begin{enumerate}
\item $[Z, K_0]=0$;
\item $[Z, P_1]=\frac{1-x}{x}[K_0, P_1]$;
\item $[Z, P_2]=\frac{c_2(y-1)}{c_1y}[K_0, P_2]$.
\end{enumerate}
\end{lemma}
\begin{proof}
If $(G/H, g)$ is a geodesic orbit space, then by Corollary \ref{GOcon}, for any $$X=X_0+X_1+X_2=(K_0, -\frac{c_2}{c_1}K_0)+(P_1, 0)+(0, P_2)\in\m,$$ there exists some $\bar{Z}=(Z, Z)\in\fh$ such that
\begin{align}
[X+\bar{Z},AX]&=[(K_0+Z+P_1, -\frac{c_2}{c_1}K_0+Z+P_2), (K_0+xP_1, -\frac{c_2}{c_1}K_0+yP_2)]\\
&=([K_0+Z+P_1, K_0+xP_1], [-\frac{c_2}{c_1}K_0+Z+P_2, -\frac{c_2}{c_1}K_0+yP_2])\in\fh.
\end{align}
Therefore, for any $K_0\in\fk, P_1\in\p_1$ and $P_2\in\p_2$, there exists some $Z\in\fk$ satisfying
\begin{equation*}
[Z, K_0]+(x-1)[K_0, P_1]+x[Z, P_1]=-\frac{c_2}{c_1}[Z, K_0]+\frac{c_2}{c_1}(1-y)[K_0, P_2]+y[Z, P_2]\in\fk,
\end{equation*}
Thus we have
\begin{enumerate}
\item $[Z, K_0]=0$;
\item $[Z, P_1]=\frac{1-x}{x}[K_0, P_1]$;
\item $[Z, P_2]=\frac{c_2(y-1)}{c_1y}[K_0, P_2]$.
\end{enumerate}
Hence, we obtain equivalent conditions in the theorem.
\end{proof}

\begin{coro}\label{coro1}
If $\frac{1-x}{x}=\frac{c_2(y-1)}{c_1y}$, then the $G$-invariant metric $g$ on $G/H$ induced by the inner product $(\ref{metric})$ is geodesic orbit.
\end{coro}
\begin{proof}
It is sufficient to take $Z=\frac{1-x}{x}K_0$ in Lemma~\ref{ec}.
\end{proof}

Next, using principal isotropy subgroups of isotropy representations, we will prove that the $G$-invariant geodesic orbit metric on $G/H$ induced by (\ref{metric}) must satisfy $\frac{1-x}{x}=\frac{c_2(y-1)}{c_1y}$.

\begin{lemma}\label{l1}
Assume that the $G$-invariant metrics on $G/H$ induced by (\ref{metric}) are geodesic orbit. If one of the principal isotropy subalgebras is either trivial or a torus in $\fk$, then $\frac{1-x}{x}=\frac{c_2(y-1)}{c_1y}$.
\end{lemma}
\begin{proof}
First, we prove the case when one of the principal isotropy subalgebras is trivial. Without loss of generality, we assume $(\fh_{\rho_1})$ is trivial. Choose $P_1$ in some $K$-orbit in general position of $\p_1$. If the $G$-invariant metric on $G/H$ induced by (\ref{metric}) is geodesic orbit, then by condition (2) in Theorem \ref{ec}, for this special $P_1\in\p_1$, $Z=\frac{1-x}{x}K_0$ because $(\fh_{\rho_1})$ is trivial. Furthermore, condition (3) in Theorem \ref{ec} implies that, for any $K_0\in\fk, P_2\in\p_2$, the following equation holds:
$$\left(\frac{c_2(y-1)}{c_1y}-\frac{1-x}{x}\right)[K_0, P_2]=0.$$
Hence $\frac{1-x}{x}=\frac{c_2(y-1)}{c_1y}$.

Then, we prove the case when one of the principal isotropy subalgebras is a torus in $\fk$. Without loss of generality, we assume $(\fh_{\rho_1})$ is a torus in $\fk$. Then we can choose $P_1$ in some $K$-orbit in general position in $\p_1$. If the $G$-invariant metric induced by (\ref{metric}) is geodesic orbit, then by condition $(2)$ in Theorem \ref{ec}, we deduce that for this special $P_1$ and for any $K_0\in\fk, P_2\in\p_2$, $Z-\frac{1-x}{x}K_0$ is in $(\fh_{\rho_1})$ which is a torus of $\fk$. Since $(\fh_{\rho_1})$ is a commutative subalgebra and every $\ad_U|_{\p_2}$ is skew symmetric for any $U\in(\fh_{\rho_1})$, there exists a basis of $\p_2$ such that $\p_2=\fl_0+\fl_1+\dots+\fl_s$, where $[(\fh_{\rho_1}), \fl_0]=0, [(\fh_{\rho_1}), \fl_i]\subset\fl_i$ and $\mathrm{dim}(\fl_i)=2$ for $i=1,\dots,s$. By condition $(3)$, for any $K_0\in\fk, P_2\in\p_2$, there exists some $Z-\frac{1-x}{x}K_0\in(\fh_{\rho_1})$ such that
$$[Z-\frac{1-x}{x}K_0, P_2]=\left(\frac{c_2(y-1)}{c_1y}-\frac{1-x}{x}\right)[K_0, P_2].$$
Hence, for any $K_0\in\fk$,
$$\left(\frac{c_2(y-1)}{c_1y}-\frac{1-x}{x}\right)[K_0, \fl_i]=[Z-\frac{1-x}{x}K_0, \fl_i]\subset\fl_i.$$
If $\frac{c_2(y-1)}{c_1y}-\frac{1-x}{x}\neq0$, then $[\fk, \fl_i]\subset\fl_i$. Since every two-dimensional representation of a compact simple Lie algebra is trivial, we have $[\fk, \p_2]=0$, which is a contradiction. As a result, we have $\frac{c_2(y-1)}{c_1y}=\frac{1-x}{x}$.
\end{proof}

%\begin{remark}\label{tori}
%Throughout the proof of the above lemma, we didn't use the property of ``maximal" of $(H_{\rho_1})$, as a result, the lemma is also true for the condition ``if one of the principal isotropy subgroups is a torus".
%\end{remark}

%\begin{remark}
%In the above proof, for the special $P_1$ on an orbit in general position of $\p_1$, we have $Z\in\ft$ such that $[Z, K_0]=0$ for any $K_0\in\fk$. In fact, we can choose some $K_0\in\fk$ such that $[K_0, T]\neq0$ for any $T\in\ft$. As a result, $Z$ must be $0$, which can also prove the above Lemma.
%\end{remark}

In summary, we have the following theorem.

\begin{theorem}\label{maintheorem}
Let $G/H=(G_1\times G_2)/\Delta K$ for $G_1\not=G_2$. Then the $G$-invariant metric on $G/H$ induced by (\ref{metric}) is geodesic orbit if and only if $\frac{1-x}{x}=\frac{c_2(y-1)}{c_1y}$. Furthermore, any geodesic orbit metric on $G/H$ is naturally reductive.
\end{theorem}

\begin{proof}
By Theorem \ref{pis}, Corollary~\ref{coro1}, Lemma \ref{l1}, we can easily obtain that any $G$-invariant metric on $G/H$ induced by (\ref{metric}) is geodesic orbit if and only if $\frac{1-x}{x}=\frac{c_2(y-1)}{c_1y}$. In the proof of Corollary~\ref{coro1}, the correspondence $X\rightarrow\bar{Z}$ can be chosen linearly and $\Ad(H)$-equivariantly. In fact, for any $X=(K_0, -\frac{c_2}{c_1}K_0)+(P_1, 0)+(0, P_2)\in\m$, there exists a linear map $\xi:\m\rightarrow\fk$ such that $\xi(X)=(\frac{1-x}{x}K_0, \frac{1-x}{x}K_0)$. By Lemma~\ref{geogra}, these geodesic orbit metrics are all naturally reductive.
\end{proof}

\section{Geodesic orbit metrics on $(G_1\times G_1)/\Delta K$}
In this case, $G/H=(G_1\times G_1)/\Delta K$, where $G_1/K$ is an isotropy irreducible space of compact type. Then $\g=\g_1\oplus\g_1$ has the following $B$-orthogonal $reductive$ $decomposition$:
\begin{equation}
\g_1\oplus\g_1=\fh\oplus\m=\fh\oplus\m_0\oplus\m_1\oplus\m_2,
\end{equation}
where $\fh=\{(Z, Z)|Z\in\fk\}$, $\m_0=\{(Z, -Z)|Z\in\fk\}$, $\m_1=\{(P, 0)|P\in\p_1\}$ and $\m_2=\{(0, P)|P\in\p_1\}$.

Note that $\m_0$ is an $\ad(\fh)$-irreducible submodule, $\m_1$ and $\m_2$ are equivalent irreducible $\ad(\fh)$-submodules. Therefore, the metric endomorphism $A$ restricted to $\m_0$ is proportional to the identity map, while restricted to $\m_1\oplus\m_2$ is a positively definite symmetric matrix under the basis of $\m_1\oplus\m_2$. Assume $A|_{\m_1\oplus\m_2}$ has two distinct eigenvalues $\lambda_1$ and $\lambda_2$, and denote by $\n_1$ and $\n_2$ the eigenspaces that belong to $\lambda_1$ and $\lambda_2$. Then $\n_1$ and $\n_2$ are also equivalent $\ad(\fh)$-irreducible submodules.

\begin{lemma}\label{lemism}
Any $\ad(\fh)$-irreducible submodule $\n$ of $\m_1\oplus\m_2$ is of the form
\begin{equation}
\n=\{(\beta P, \gamma \theta(P)|P\in\p_1, \beta^2+\gamma^2\neq0\},
\end{equation}
for some fixed $\beta,\gamma\in\mathbb{R}$, where $\theta$ is an $\ad(\fh)$-equivariant isomorphism on $\p_1$.
\end{lemma}
\begin{proof}
Let $\n$ be an irreducible $\ad(\fh)$-submodule in $\m_1\oplus\m_2$. Consider the linear maps $f_i$: $\m_1\oplus\m_2\rightarrow\m_i(i=1,2)$ defined by $f_i(P)=P_i$, where $P=(P_1, P_2)$ with $P_i\in\m_i(i=1,2)$. Since $\n, \m_1$ and $\m_2$ are $\ad(\fh)$-invariant, for any $Z\in\fh$, we have $f_i([Z, P])=[Z, f_i(P)](i=1,2)$. Since $\m_i(i=1,2)$ are $\ad(\fh)$-irreducible, this implies $f_i(\n)$ is either $0$ or $\m_i$ for $i=1,2$. If $f_i(\n)=0$, then $f_j(\n)=\m_j$ for $\{i, j\}=\{1, 2\}$, the lemma is obviously true. Otherwise, $f_i(i=1, 2)$ are bijective. Consider the linear mapping $\theta=f_{2}\circ f_{1}^{-1}:\p_1\rightarrow\p_1$, it is easy to verified that $\theta$ is $\ad(\fh)$-equivariant. That is $P_2=\theta(P_1)$. This completes the proof of this lemma.
\end{proof}

In \cite{AlNi}, the authors found a very useful property of the metric endomorphism $A$ determined by a geodesic orbit metric on homogeneous space $G/H$. We will recall their results along with the proof in the following.
\begin{prop}[Corollary 3, \cite{AlNi}]\label{prop1}
Let $(G/H, g)$ be a compact geodesic orbit space with the reductive decomposition (\ref{rd}) and
$$\m=\ffb_1+\ffb_2+\cdots+\ffb_k$$
the $A$-eigenspace decomposition such that $A|_{\ffb_i}=\lambda_i1_{\ffb_i}$. Then for any $\ad(\fh)$-submodules $\q_i\subset\ffb_i$, $\q_j\subset\ffb_j, i\neq j$, we have
$$[\q_i, \q_j]\subset\q_i+\q_j.$$
\end{prop}
\begin{proof}
Let $X_i\in\q_i, X_j\in\q_j$ be eigenvectors of $A$ with different eigenvalues $\lambda_i,\lambda_j(i\neq j)$. By Lemma~\ref{lemma=0}, there exists $W\in\fh$ such that
$$[W+X_i+X_j, A(X_i+X_j)]=0,$$
that is,
$$(\lambda_i-\lambda_j)[X_i, X_j]=[W, \lambda_iX_i+\lambda_jX_j].$$
Since $\q_i$ and $\q_j$ are $\ad(\fh)$-submodules, the right hand of above identity belongs to $\q_i+\q_j$. This completes the proof of this proposition.
\end{proof}

Let $\m=\m_0\oplus\n_1\oplus\n_2$ be any decomposition of $\m$ into $\ad(\Delta(\fk))$-irreducible submodules. If $\n_i\neq\m_i(i=1,2)$, then by Lemma~\ref{lemism}, $\n_1=\{(\beta_1 P, \gamma_1 \theta_1(P)|P\in\p_1, \beta_1^2+\gamma_1^2\neq0\}$ and $\n_2=\{(\beta_2 P, \gamma_2 \theta_2(P)|P\in\p_1, \beta_2^2+\gamma_2^2\neq0\}$, where $\beta_1\beta_2\gamma_1\gamma_2\neq0$. By Proposition \ref{prop1}, for any $P_1, P_2\in\p_1$,
\begin{equation}
[(\beta_1 P_1, \gamma_1 \theta_1(P_1)), (\beta_2 P_2, \gamma_2 \theta_2(P_2))]=(\beta_1\beta_2[P_1, P_2], \gamma_1\gamma_2[\theta_1(P_1), \theta_2(P_2)])\in\n_1\oplus\n_2\cong\p_1\oplus\p_1,
\end{equation}
which implies that $[\p_1, \p_1]\subset\p_1$. Hence, $\p_1$ is an ideal in $\g_1$, which is a contradiction. Therefore, if $\n_i\neq\m_i, (i=1,2)$, then $A|_{\m_1\oplus\m_2}$ is proportional to the identity map. Along with the case when $\n_i=\m_i, (i=1,2)$, we have

%From the discussion above, we obtain that if $A$ induces a geodesic orbit metric, then $A|_{\m_1\oplus\m_2}$ is proportional to the identity map. Therefore, in this case, any $G$-invariant geodesic orbit metric on $G/H$ is of the following form
%\begin{equation}\label{metric1113}
%(\ ,\ )=B|_{\m_0}+xB|_{\m_1}+xB|_{\m_2},\quad x\in\mathbb{R}^+.
%\end{equation}
%Therefore, we have:
\begin{theorem} Any $G$-invariant geodesic orbit metric on $G/H=(G_1\times G_1)/\Delta K$ is diagonal with respect to the decomposition $\m=\m_0\oplus\m_1\oplus\m_2$, that is,
\begin{equation}\label{metric1112}
(\ ,\ )=B|_{\m_0}+xB|_{\m_1}+yB|_{\m_2},\quad x,y\in\mathbb{R}^+.
\end{equation}
\end{theorem}
Let $(G/H, g)$ be a geodesic orbit space. We can choose $K_0\in\fk$ and $P_1=P_2=P\in\p_1$  satisfying $[K_0, P]\neq0$. Similar to Lemma~\ref{ec}, there exists some $Z\in\fk$ such that
\begin{enumerate}
\item $[Z, K_0]=0$;
\item $[Z, P]=\frac{1-x}{x}[K_0, P]$;
\item $[Z, P]=\frac{y-1}{y}[K_0, P]$.
\end{enumerate}
Therefore, $\frac{1-x}{x}=\frac{y-1}{y}.$ Conversely, it's easy to verify that the metric induced by the inner product of the form (\ref{metric1112}) with $\frac{1-x}{x}=\frac{y-1}{y}$ is geodesic orbit. In fact, we can choose $Z=\frac{y-1}{y}K_0$. Thus we have,
\begin{theorem}\label{maintheorem1}
Let the notation be as above. Any geodesic orbit metric on $G/H=(G_1\times G_1)/\Delta K$ is of the form (\ref{metric1112}) with $\frac{1-x}{x}=\frac{y-1}{y}$. Furthermore, any geodesic orbit metric is naturally reductive.
\end{theorem}

In summary, Theorem~\ref{mainthm} follows from Theorems~\ref{maintheorem} and~\ref{maintheorem1}.

%\begin{remark}
%In Cartesian coordinate system, the point $(x, y)$ is on one branch of the hyperbolic $y=\frac{x}{2x-1}$(see the picture below). By Lemma \ref{nr}, the geodesic orbit metric is naturally reductive if and only if $x=y=1$.
%\end{remark}
%\begin{center}
%\includegraphics[width=14em]{graph1.jpeg}
%\end{center}

%\begin{remark}
%The geodesic orbit metric is determined by the pair $(x, y)$ in (\ref{metric}) satisfying $\frac{1-x}{x}=\frac{c_2(y-1)}{\alpha c_1y}$, $x,y>0$, i.e. $y=\frac{c_2x}{(\alpha c_1+c_2)x-\alpha c_1}$, $x,y>0$, which denotes exactly the branch of a hyperbola passing through the point $(1,1)$ in the $(x, y)$-plane. For example, see the picture below for $c_1=c_2$, which is the case for $G_1=G_2$.
%\end{remark}
%\begin{center}
%\includegraphics[width=14em]{graph1.jpeg}
%\end{center}

%\begin{remark}
%Notice that all the homogeneous spaces $G/H$ constructed here are not weakly symmetric spaces with reductive isometry groups (see \cite{Ya}). As a result, we obtain new non-naturally reductive geodesic orbit metrics on homogeneous spaces with semisimple (non-simple) isometry groups which are not weakly symmetric spaces.
%\end{remark}

\section{Acknowledgements} This work was partially supported by National Natural Science Foundation of China (11571182, 11901300 and 11931009), Natural Science Foundation of Tianjin (19JCYBJC30600) and Natural Science Research of Jiangsu Education Institutions of China (No. 19KJB110015). The first author would like to thank Prof. Yu.G. Nikonorov for very useful discussions and suggestions. Also he would like to thank Prof. Yu.G. Nikonorov and Y. Nikolayevsky for useful discussions during their visit at Chern Institute of Mathematics via Chern Institute of Mathematics visiting scholar programme.

\section*{Appendix I}
The classification of $G_1\times G_2/\Delta(K)$ for $G_1\not=G_2$ is given as follows. Here $\rho_i$ $(i=1,2)$ denote the isotropy representations of $G_i/K$. Let $\fk$ be the Lie algebra of $K$ and let $\{\alpha_1, \dots, \alpha_r\}$ be a system of simple roots of $\fk$. Identify the isotropy representation with its complexification and denote by $\varphi_1, \dots, \varphi_r$ the fundamental  representations corresponding to $\alpha_1, \dots, \alpha_r$ respectively.

\textbf{Type A:}
$\fk=\mathfrak{su}(n)$.

A.1 \quad  $\Delta\mathrm{SU}(n)/\mathbb{Z}_n< \mathrm{Spin}(n^{2}-1)\times \mathrm{SU}(\frac{n(n-1)}{2})/\mathbb{Z}_n$, $n\geq 5, n \mbox{ odd}.$ Here
$$\rho_1=(\varphi_1+2\varphi_{n-1})\oplus(2\varphi_1+\varphi_{n-1}),\ \ \ \rho_2=\varphi_2+\varphi_{n-2}.$$

A.2 \quad  $\Delta\mathrm{SU}(n)/\mathbb{Z}_n< \mathrm{SO}(n^{2}-1)\times \mathrm{SU}(\frac{n(n-1)}{2})/\mathbb{Z}_{\frac{n}{2}}$, $n\geq 6, n \mbox{ even}.$ Here
$$\rho_1=(\varphi_1+2\varphi_{n-1})\oplus(2\varphi_1+\varphi_{n-1}),\ \ \ \rho_2=\varphi_2+\varphi_{n-2}.$$

A.3 \quad  $\Delta\mathrm{SU}(n)/\mathbb{Z}_n< \mathrm{Spin}(n^{2}-1)\times \mathrm{SU}(\frac{n(n+1)}{2})/\mathbb{Z}_n$, $n\geq 3, n \mbox{ odd}.$ Here
\begin{equation*}
\rho_1=\left\{
             \begin{array}{ll}
             3\varphi_2\oplus3\varphi_1 & n=3,  \\
              (\varphi_2+2\varphi_{n})\oplus(2\varphi_1+\varphi_{n-1})& n>3\ \mbox{and odd},
             \end{array}
\right. \rho_2=2\varphi_1+2\varphi_{n-1}.
\end{equation*}

A.4 \quad  $\Delta\mathrm{SU}(n)/\mathbb{Z}_n< \mathrm{SO}(n^{2}-1)\times \mathrm{SU}(\frac{n(n+1)}{2})/\mathbb{Z}_{\frac{n}{2}}$, $n\geq 4, n \mbox{ even}.$ Here
$$\rho_1=(\varphi_1+2\varphi_{n-1})\oplus(2\varphi_1+\varphi_{n-1}),\ \ \ \rho_2=2\varphi_1+2\varphi_{n-1}.$$

A.5 \quad  $\Delta\mathrm{SU}(n)/\mathbb{Z}_n< \mathrm{SU}(\frac{n(n-1)}{2})/\mathbb{Z}_m \times \mathrm{SU}(\frac{n(n+1)}{2})/\mathbb{Z}_m$, $ n\geq 5$, $m=\left\{\begin{array}{ll}\frac{n}{2}, &n\mbox{ even},\\ n,&n \mbox{ odd}.\end{array}\right.$ Here
$$\rho_1=\varphi_2+\varphi_{n-2},\ \ \ \rho_2=2\varphi_1+2\varphi_{n-1}.$$

A.6  \quad  $\Delta\mathrm{SU}(3)/\mathbb{Z}_3< \mathrm{E}_6\times \mathrm{SU}(6)/\mathrm{Z}_3.$ Here $\rho_1=(\varphi_1+4\varphi_2)\oplus(4\varphi_1+\varphi_2)$, $\rho_2=2\varphi_1+2\varphi_2$.

A.7  \quad  $\Delta\mathrm{SU}(3)/\mathbb{Z}_3< \mathrm{E}_6\times \mathrm{Spin}(8).$ Here $\rho_1=(\varphi_1+4\varphi_2)\oplus(4\varphi_1+\varphi_2)$, $\rho_2=3\varphi_2\oplus3\varphi_1$.

A.8  \quad  $\Delta\mathrm{SU}(3)/\mathbb{Z}_3< \mathrm{E}_6\times \mathrm{E}_7.$ Here $\rho_1=(\varphi_1+4\varphi_2)\oplus(4\varphi_1+\varphi_2)$, $\rho_2=4\varphi_1+4\varphi_2$.

A.9  \quad  $\Delta\mathrm{SU}(3)/\mathbb{Z}_3< \mathrm{SU}(6)/\mathbb{Z}_3\times \mathrm{Spin}(8).$ Here $\rho_1=2\varphi_1+2\varphi_2$, $\rho_2=3\varphi_1+3\varphi_2$.

A.10  \quad  $\Delta\mathrm{SU}(3)/\mathbb{Z}_3< \mathrm{SU}(6)/\mathbb{Z}_3\times \mathrm{E}_7.$ Here $\rho_1=2\varphi_1+2\varphi_2$, $\rho_2=4\varphi_1+4\varphi_2$.

A.11  \quad  $\Delta\mathrm{SU}(3)/\mathbb{Z}_3< \mathrm{Spin}(8)\times \mathrm{E}_7.$ Here $\rho_1=3\varphi_2\oplus3\varphi_1$, $\rho_2=4\varphi_1+4\varphi_2$.

A.12  \quad  $\Delta\mathrm{SU}(4)/\mathbb{Z}_4< \mathrm{SO}(20)\times \mathrm{SO}(15).$ Here $\rho_1=\varphi_1+2\varphi_2+\varphi_3$, $\rho_2=(\varphi_1+2\varphi_2)\oplus(2\varphi_1+\varphi_2)$.

A.13  \quad  $\Delta\mathrm{SU}(4)/\mathbb{Z}_4< \mathrm{SO}(20)\times \mathrm{SU}(10)/\zz.$ Here $\rho_1=\varphi_1+2\varphi_2+\varphi_3$, $\rho_2=2\varphi_1+2\varphi_3$.

A.14  \quad  $\Delta\mathrm{SU}(4)/\mathbb{Z}_4< \mathrm{SO}(15)\times \mathrm{SU}(10)/\zz.$ Here $\rho_1=(\varphi_1+2\varphi_2)\oplus(2\varphi_1+\varphi_2)$, $\rho_2=2\varphi_1+2\varphi_3$.

A.15  \quad  $\Delta\mathrm{SU}(6)/\mathbb{Z}_6< \mathrm{SO}(35)\times \mathrm{Sp}(10)/\zz.$ Here $\rho_1=2\varphi_1+2\varphi_5$, $\rho_2=2\varphi_3$.

A.16  \quad  $\Delta\mathrm{SU}(6)/\mathbb{Z}_6< \mathrm{SO}(35)\times \mathrm{SU}(15)/\mathbb{Z}_3.$ Here $\rho_1=2\varphi_1+2\varphi_5$, $\rho_2=\varphi_2+\varphi_4$.

A.17  \quad  $\Delta\mathrm{SU}(6)/\mathbb{Z}_6< \mathrm{SO}(35)\times \mathrm{SU}(21)/\mathbb{Z}_3.$ Here $\rho_1=2\varphi_1+2\varphi_5$, $\rho_2=2\varphi_1+2\varphi_5$.

A.18  \quad  $\Delta\mathrm{SU}(6)/\mathbb{Z}_6< \mathrm{Sp}(10)/\zz \times \mathrm{SU}(15)/\mathbb{Z}_3.$ Here $\rho_1=2\varphi_3$, $\rho_2=\varphi_2+\varphi_4$.

A.19  \quad  $\Delta\mathrm{SU}(6)/\mathbb{Z}_6< \mathrm{Sp}(10)/\zz \times \mathrm{SU}(21)/\mathbb{Z}_3.$ Here $\rho_1=2\varphi_3$, $\rho_2=2\varphi_1+2\varphi_5$.

A.20  \quad  $\Delta\mathrm{SU}(6)/\mathbb{Z}_6< \mathrm{SU}(21)/\mathbb{Z}_3 \times \mathrm{SU}(15)/\mathbb{Z}_3.$ Here $\rho_1=2\varphi_1+2\varphi_5$, $\rho_2=\varphi_2+\varphi_4$.

A.21  \quad  $\Delta\mathrm{SU}(8)/\mathbb{Z}_8< \mathrm{SO}(70)/\mathbb{Z}_2 \times \mathrm{SO}(63).$ Here $\rho_1=\varphi_3+\varphi_5$, $\rho_2=(\varphi_2+2\varphi_7)\oplus(2\varphi_1+\varphi_6)$.

A.22  \quad  $\Delta\mathrm{SU}(8)/\mathbb{Z}_8< \mathrm{SO}(70)/\mathbb{Z}_2 \times \mathrm{SU}(28)/\mathbb{Z}_4.$ Here $\rho_1=\varphi_3+\varphi_5$, $\rho_2=\varphi_2+\varphi_6$.

A.23  \quad  $\Delta\mathrm{SU}(8)/\mathbb{Z}_8< \mathrm{SO}(70)/\mathbb{Z}_2 \times \mathrm{SU}(36)/\mathbb{Z}_4.$ Here $\rho_1=\varphi_3+\varphi_5$, $\rho_2=2\varphi_1+2\varphi_7$.

A.24  \quad  $\Delta\mathrm{SU}(8)/\mathbb{Z}_8< \mathrm{SO}(63) \times \mathrm{SU}(28)/\mathbb{Z}_4.$ Here $\rho_1=(\varphi_2+2\varphi_7)\oplus(2\varphi_1+\varphi_6)$, $\rho_2=\varphi_2+\varphi_6$.

A.25  \quad  $\Delta\mathrm{SU}(8)/\mathbb{Z}_8< \mathrm{SO}(63) \times \mathrm{SU}(36)/\mathbb{Z}_4.$ Here $\rho_1=(\varphi_2+2\varphi_7)\oplus(2\varphi_1+\varphi_6)$, $\rho_2=2\varphi_1+2\varphi_7$.

A.26  \quad  $\Delta\mathrm{SU}(8)/\mathbb{Z}_8< \mathrm{SU}(28)/\mathbb{Z}_4 \times \mathrm{SU}(36)/\mathbb{Z}_4.$ Here $\rho_1=\varphi_2+\varphi_6$, $\rho_2=2\varphi_1+2\varphi_7$.

A.27  \quad  $\Delta\mathrm{SU}(9)/\mathbb{Z}_9< \mathrm{Spin}(80) \times \mathrm{SU}(36)/\mathbb{Z}_9.$ Here $\rho_1=(\varphi_2+2\varphi_8)\oplus(2\varphi_1+\varphi_8)$, $\rho_2=\varphi_2+\varphi_7$.

A.28  \quad  $\Delta\mathrm{SU}(9)/\mathbb{Z}_9< \mathrm{Spin}(80) \times \mathrm{SU}(45)/\mathbb{Z}_9.$ Here $\rho_1=(\varphi_2+2\varphi_8)\oplus(2\varphi_1+\varphi_8)$, $\rho_2=2\varphi_1+2\varphi_8$.

A.29  \quad  $\Delta\mathrm{SU}(9)/\mathbb{Z}_9< \mathrm{SU}(36)/\mathbb{Z}_9 \times \mathrm{SU}(45)/\mathbb{Z}_9.$ Here $\rho_1=\varphi_2+\varphi_7$, $\rho_2=2\varphi_1+2\varphi_8$.

\textbf{Type B:}
$\fk=\mathfrak{so}(n)$.

B.1 \quad $\Delta \mathrm{SO}(n)< \mathrm{SO}(n+1)\times \mathrm{SO}(\frac{n(n-1)}{2})$, $ n\geq 7.$ Here
$$\rho_1=\varphi_1,\ \ \
\rho_2=\left\{
             \begin{array}{ll}
             \varphi_1+2\varphi_3,&n=7, \\
             \varphi_1+\varphi_3+\varphi_4, &n=8,\\
              \varphi_1+\varphi_3,& n>8.\end{array}\right.$$

B.2 \quad $\Delta \mathrm{SO}(n)< \mathrm{SO}(n+1)\times \mathrm{SO}(\frac{(n-1)(n+2)}{2})$, $  n\geq 5.$ Here
$$\rho_1=\varphi_1,\ \ \
\rho_2=\left\{
             \begin{array}{ll}
              2\varphi_1+2\varphi_3, &n=5,\\
              \varphi_1+2\varphi_2+\varphi_3, &n=6\ (D_3=A_3),\\
             2\varphi_1+\varphi_2, &n>6.\end{array}\right.$$

B.3 \quad $\Delta \mathrm{SO}(n)< \mathrm{SO}(n+1)\times \mathrm{SU}(n)$, $  n\geq 3.$ Here
$$\rho_1=\varphi_1,\ \ \
\rho_2=\left\{
             \begin{array}{ll}
              4\varphi_1, &n=3, \\
              2\varphi_1\otimes2\varphi_1,&n=4\ (D_2=A_1\oplus A_1),  \\
             2\varphi_1,&n>4.
             \end{array}\right.
$$

B.4 \quad $\Delta \mathrm{SO}(n)<  \mathrm{SO}(\frac{n(n-1)}{2}) \times \mathrm{SO}(\frac{(n-1)(n+2)}{2})$, $  n\geq 7.$ Here
$$\rho_1=\left\{
             \begin{array}{ll}
             \varphi_1+2\varphi_3, &n=7,  \\
   \varphi_1+\varphi_3+\varphi_4, &n=8,  \\
              \varphi_1+\varphi_3,&n>8,\end{array}\right.
\ \ \  \rho_2=2\varphi_1+\varphi_2.$$

B.5 \quad $\Delta \mathrm{SO}(n)<  \mathrm{SO}(\frac{n(n-1)}{2}) \times \mathrm{SU}(n)$, $  n\geq 7.$ Here
$$\rho_1=\left\{
             \begin{array}{ll}
             \varphi_1+2\varphi_3,& n=7, \\
             \varphi_1+\varphi_3+\varphi_4,&n=8, \\
              \varphi_1+\varphi_3,&n>8,\end{array}\right. \ \ \ \rho_2=2\varphi_1.$$

B.6 \quad $\Delta \mathrm{SO}(n)<  \mathrm{SO}(\frac{(n-1)(n+2)}{2}) \times \mathrm{SU}(n)$, $  n\geq 5.$ Here
$$\rho_1=\left\{
             \begin{array}{ll}
              2\varphi_1+2\varphi_3,&n=5, \\
              \varphi_1+2\varphi_2+\varphi_3,&n=6\ (D_3=A_3), \\
             2\varphi_1+\varphi_2,&n>6,\end{array}\right.\ \ \ \rho_2=2\varphi_1.
$$

B.7 \quad $\Delta \mathrm{SO}(3)<  \mathrm{SU}(3) \times \mathrm{Sp}(2)/\mathbb{Z}_2.$ Here $\rho_1=4\varphi_2$, $\rho_2=6\varphi_1$.

B.8 \quad $\Delta \mathrm{SO}(3)<  \mathrm{SU}(3) \times \mathrm{G}_{2}.$ Here $\rho_1=4\varphi_2$, $\rho_2=10\varphi_1$.

B.9 \quad $\Delta \mathrm{SO}(3)<  \mathrm{Sp}(2)/\zz \times \mathrm{G}_{2}.$ Here $\rho_1=6\varphi_1$, $\rho_2=10\varphi_1$.

B.10 \quad $\Delta \mathrm{SO}(9)<  \mathrm{SO}(10) \times \mathrm{Spin}(16)/\zz.$ Here $\rho_1=\varphi_1$, $\rho_2=\varphi_3$.

B.11  \quad $\Delta \mathrm{SO}(9)<  \mathrm{SO}(10) \times \mathrm{SO}(36).$ Here $\rho_1=\varphi_1$, $\rho_2=\varphi_1+\varphi_3$.

B.12  \quad $\Delta \mathrm{SO}(9)<  \mathrm{SO}(10) \times \mathrm{SO}(44).$ Here $\rho_1=\varphi_1$, $\rho_2=2\varphi_1+\varphi_2$.

B.13  \quad $\Delta \mathrm{SO}(9)<  \mathrm{SO}(10) \times \mathrm{SU}(9).$ Here $\rho_1=\varphi_1$, $\rho_2=2\varphi_1$.

B.14  \quad $\Delta \mathrm{SO}(9)<  \mathrm{Spin}(16)/\zz \times \mathrm{SO}(36).$ Here $\rho_1=\varphi_3$, $\rho_2=\varphi_1+\varphi_3$.

B.15  \quad $\Delta \mathrm{SO}(9)<  \mathrm{Spin}(16)/\zz \times \mathrm{SO}(44).$ Here $\rho_1=\varphi_3$, $\rho_2=2\varphi_1+\varphi_2$.

B.16  \quad $\Delta \mathrm{SO}(9)<  \mathrm{Spin}(16)/\zz \times \mathrm{SU}(9).$ Here $\rho_1=\varphi_3$, $\rho_2=2\varphi_1$.

B.17  \quad $\Delta \mathrm{SO}(9)<  \mathrm{SO}(36) \times \mathrm{SO}(44).$ Here $\rho_1=\varphi_1+\varphi_3$, $\rho_2=2\varphi_1+\varphi_2$.

B.18  \quad $\Delta \mathrm{SO}(9)<  \mathrm{SO}(36) \times \mathrm{SU}(9).$ Here $\rho_1=\varphi_1+\varphi_3$, $\rho_2=2\varphi_1$.

B.19  \quad $\Delta \mathrm{SO}(9)<  \mathrm{SO}(44) \times \mathrm{SU}(9).$ Here $\rho_1=2\varphi_1+\varphi_2$, $\rho_2=2\varphi_1$.

B.20  \quad $\Delta \mathrm{SO}(10)<  \mathrm{SO}(11) \times \mathrm{SU}(10).$ Here $\rho_1=\varphi_1$, $\rho_2=2\varphi_1$.

B.21  \quad $\Delta \mathrm{SO}(10)/\zz<  \mathrm{SU}(16)/\zzzz \times \mathrm{SO}(45).$ Here $\rho_1=\varphi_4+\varphi_5$, $\rho_2=\varphi_1+\varphi_3$.

B.22  \quad $\Delta \mathrm{SO}(10)/\zz<  \mathrm{SU}(16)/\zzzz \times \mathrm{SO}(54).$ Here $\rho_1=\varphi_4+\varphi_5$, $\rho_2=2\varphi_1+\varphi_2$.

B.23  \quad $\Delta \mathrm{SO}(10)/\zz<  \mathrm{SO}(54) \times \mathrm{SO}(45).$ Here $\rho_1=2\varphi_1+\varphi_2$, $\rho_2=\varphi_1+\varphi_3$.

B.24  \quad $\Delta \mathrm{SO}(12)<  \mathrm{SU}(12) \times \mathrm{SO}(13).$ Here $\rho_1=2\varphi_1$, $\rho_2=\varphi_1$.

B.25  \quad $\Delta \mathrm{SO}(12)/\zz<  \mathrm{Sp}(16)/\zz \times \mathrm{Spin}(66).$ Here $\varphi_1=2\varphi_5$, $\rho_2=\varphi_1+\varphi_3$.

B.26  \quad $\Delta \mathrm{SO}(12)/\zz<  \mathrm{Sp}(16)/\zz \times \mathrm{Spin}(77).$ Here $\varphi_1=2\varphi_5$, $\rho_2=2\varphi_1+\varphi_2$.

B.27  \quad $\Delta \mathrm{SO}(12)/\zz<  \mathrm{Spin}(77) \times \mathrm{Spin}(66).$ Here $\rho_1=2\varphi_1+\varphi_2$, $\rho_2=\varphi_1+\varphi_3$.

B.28  \quad $\Delta \mathrm{SO}(16)<  \mathrm{SO}(17) \times \mathrm{E}_8.$ Here $\rho_1=\varphi_1$, $\rho_2=\varphi_7$.

B.29  \quad $\Delta \mathrm{SO}(16)<  \mathrm{SO}(17) \times \mathrm{SU}(16).$ Here $\rho_1=\varphi_1$, $\rho_2=2\varphi_1$.

B.30  \quad $\Delta \mathrm{SO}(16)<  \mathrm{SU}(16) \times \mathrm{E}_8.$ Here $\rho_1=2\varphi_1$, $\rho_2=\varphi_7$.

B.31  \quad $\Delta \mathrm{SO}(16)/\zz<  \mathrm{Spin}(128)/\zz \times \mathrm{Spin}(120).$ Here $\rho_1=\varphi_6$, $\rho_2=\varphi_1+\varphi_3$.

B.32  \quad $\Delta \mathrm{SO}(16)/\zz<  \mathrm{Spin}(128)/\zz \times \mathrm{Spin}(135).$ Here $\rho_1=\varphi_6$, $\rho_2=2\varphi_1+\varphi_2$.

B.33  \quad $\Delta \mathrm{SO}(16)/\zz<  \mathrm{Spin}(135) \times \mathrm{Spin}(120).$ Here $\rho_1=2\varphi_1+\varphi_2$, $\rho_2=\varphi_1+\varphi_3$.

\textbf{Type C:}
$\fk=\mathfrak{sp}(n)$.

C.1 \quad  $\Delta\mathrm{Sp}(n)/\zz < \mathrm{Spin}((n-1)(2n+1))\times \mathrm{Spin}(n(2n+1))$, $n\geq 3,\ n\equiv0\ mod\ 4.$ Here $\rho_1=\varphi_1+\varphi_3$, $\rho_2=2\varphi_1+\varphi_2$.

C.2 \quad  $\Delta\mathrm{Sp}(n)/\zz < \mathrm{Spin}((n-1)(2n+1))\times \mathrm{SO}(n(2n+1))$, $n\geq 3,\ n\equiv1\ mod\ 4.$ Here $\rho_1=\varphi_1+\varphi_3$, $\rho_2=2\varphi_1+\varphi_2$.

C.3 \quad  $\Delta\mathrm{Sp}(n)/\zz < \mathrm{SO}((n-1)(2n+1))\times \mathrm{SO}(n(2n+1))$, $n\geq 3,\ n\equiv2\ mod\ 4.$ Here $\rho_1=\varphi_1+\varphi_3$, $\rho_2=2\varphi_1+\varphi_2$.

C.4 \quad  $\Delta\mathrm{Sp}(n)/\zz < \mathrm{SO}((n-1)(2n+1))\times \mathrm{Spin}(n(2n+1))$, $n\geq 3,\ n\equiv3\ mod\ 4.$ Here $\rho_1=\varphi_1+\varphi_3$, $\rho_2=2\varphi_1+\varphi_2$.

C.5 \quad  $\Delta\mathrm{Sp}(3)/\zz < \mathrm{Sp}(7)/\zz \times \mathrm{SO}(14).$ Here $\rho_1=2\varphi_3$, $\rho_2=\varphi_1+\varphi_3$.

C.6 \quad  $\Delta\mathrm{Sp}(3)/\zz < \mathrm{Sp}(7)/\zz \times \mathrm{Spin}(21).$ Here $\rho_1=2\varphi_3$, $\rho_2=2\varphi_1+\varphi_2$.

C.7 \quad  $\Delta\mathrm{Sp}(3)/\zz < \mathrm{SO}(14) \times \mathrm{Spin}(21).$ Here $\rho_1=\varphi_1+\varphi_3$, $\rho_2=2\varphi_1+\varphi_2$.

C.8 \quad  $\Delta\mathrm{Sp}(4)/\zz < \mathrm{Spin}(27) \times \mathrm{Spin}(36).$ Here $\rho_1=\varphi_1+\varphi_3$, $\rho_2=2\varphi_1+\varphi_2$.

C.9 \quad  $\Delta\mathrm{Sp}(4)/\zz < \mathrm{Spin}(27) \times \mathrm{Spin}(42).$ Here $\rho_1=\varphi_1+\varphi_3$, $\rho_2=2\varphi_3$.

C.10 \quad  $\Delta\mathrm{Sp}(4)/\zz < \mathrm{Spin}(27) \times \mathrm{E}_6.$ Here $\rho_1=\varphi_1+\varphi_3$, $\rho_2=\varphi_4$.

C.11 \quad  $\Delta\mathrm{Sp}(4)/\zz < \mathrm{Spin}(36) \times \mathrm{Spin}(42).$ Here $\rho_1=2\varphi_1+\varphi_2$, $\rho_2=2\varphi_3$.

C.12 \quad  $\Delta\mathrm{Sp}(4)/\zz < \mathrm{Spin}(36) \times \mathrm{E}_6.$ Here $\rho_1=2\varphi_1+\varphi_2$, $\rho_2=\varphi_4$.

C.13 \quad  $\Delta\mathrm{Sp}(4)/\zz < \mathrm{Spin}(42) \times \mathrm{E}_6.$ Here $\rho_1=2\varphi_3$, $\rho_2=\varphi_4$.

\textbf{Type D:}
$\fk$ is an exceptional Lie algebra.

D.1 \quad  $\Delta \mathrm{G}_{2}< \mathrm{Spin}(7)\times \mathrm{Spin}(14).$ Here $\rho_1=\varphi_1$, $\rho_2=3\varphi_1$.

D.2 \quad  $\Delta \mathrm{G}_{2}< \mathrm{Spin}(7)\times \mathrm{E}_{6}.$ Here $\rho_1=\varphi_1$, $\rho_2=\varphi_1+\varphi_2$.

D.3 \quad  $\Delta \mathrm{G}_{2}< \mathrm{Spin}(14)\times \mathrm{E}_{6}.$ Here $\rho_1=3\varphi_1$, $\rho_2=\varphi_1+\varphi_2$.

D.4 \quad  $\Delta \mathrm{F}_{4}< \mathrm{Spin}(26)\times \mathrm{Spin}(52).$ Here $\rho_1=\varphi_2$, $\rho_2=\varphi_3$.

D.5 \quad  $\Delta \mathrm{F}_{4}< \mathrm{Spin}(26)\times \mathrm{E}_{6}.$ Here $\rho_1=\varphi_2$, $\rho_2=\varphi_1$.

D.6 \quad  $\Delta \mathrm{F}_{4}< \mathrm{Spin}(52)\times \mathrm{E}_{6}.$ Here $\rho_1=\varphi_3$, $\rho_2=\varphi_1$.

D.7 \quad  $\Delta \mathrm{E}_{6}/\mathbb{Z}_3 < \mathrm{Spin}(78)\times \mathrm{SU}(27)/\mathbb{Z}_3.$ Here $\rho_1=\varphi_3$, $\rho_2=\varphi_1+\varphi_5$.

D.8 \quad  $\Delta \mathrm{E}_{7}/\zz < \mathrm{SO}(133)\times \mathrm{Sp}(28)/\zz.$ Here $\rho_1=\varphi_5$, $\rho_2=2\varphi_6$.

\newpage
\section*{Appendix II}
\begin{theorem}[\cite{HiHi} Table A]\label{HiHi}
All real irreducible representations of simple Lie groups with non-trivial principal isotropy subgroups are listed in the following tables. The notation is the same as that in the previous section.
\end{theorem}

\begin{footnotesize}
\quad I. $\fk=\begin{picture}(6,6)
\put(0,0){$A_r:$}
\put(20,0){$\circ$}
\put(23.5,2){\line(1,0){20}}
\put(43.5,0){$\circ$}
\put(48,2){\line(1,0){4}}
\put(54,1.5){$\ldots$}
%\put(50,1.5){$\ldots$}
\put(67,2){\line(1,0){6}}
\put(73,0){$\circ$}
\put(77,2){\line(1,0){20}}
%\put(65.3,1.3){\line(1,0){10}}
\put(96,0){$\circ$}
\put(19,-6){$\alpha_1$}
\put(41.5,-6){$\alpha_2$}
\put(69,-6){$\alpha_{r-1}$}
\put(94,-6){$\alpha_r$}
\end{picture}
\quad\\$
\begin{center}
\begin{tabularx}{\textwidth}{|c|c|c|Y|}
\hline
rank $r$ &\hspace{5mm}The representation $\rho$ \hspace{5mm} &\hspace{10mm} The group $K$ \hspace{10mm} & $(H_{\rho})$\\
\hline
\multirow{2}*{$r=1$} & $S^2\varphi_1$ & $\SO(3)$ & maximal tori\\
\cline{2-4}
& $S^4\varphi_1$ & $\SO(3)$ & maximal $\mathbb{Z}_2$-tori$\approx\mathbb{Z}_2\oplus\mathbb{Z}_2$\\
\hline
\multirow{3}*{$r\geq2$} & $\varphi_1+\varphi_r$ & $\SU(r+1)$ & $\SU(r)$ \\
\cline{2-4}
& $S^2\varphi_1+S^2\varphi_r$ & $\SU(r+1)$ & maximal $\mathbb{Z}_2$-tori\\
\cline{2-4}
& $\varphi_1\otimes\varphi_r$ & $\SU(r+1)$ & maximal tori\\
\hline
\multirow{2}*{$r=3$} & $\varphi_2$ & $\SO(6)$ & $\SO(5)$\\
\cline{2-4}
& $S^2\varphi_2-\theta$ & $\SO(6)$ & maximal $\mathbb{Z}_2$-tori\\
\hline
\multirow{1}*{$r\geq4$ even}& \multirow{1}*{$\varphi_2+\varphi_{r-1}$} & \multirow{1}*{$\SU(r+1)$} & $\SU(2)\times\dots\times\SU(2)$\quad  ($\big[\frac{r+1}{2}\big]$ copies) \\
\hline
\multirow{1}*{$r>4$ odd} & \multirow{1}*{$\varphi_2+\varphi_{r-1}$} & \multirow{1}*{$\SU(r+1)/\mathbb{Z}_2$} & $(\SU(2)\times\dots\times\SU(2))/\mathbb{Z}_2$\quad($\frac{r+1}{2}$ copies) \\
\hline
\multirow{1}*{$r=5$} & \multirow{1}*{$2\varphi_3$} & \multirow{1}*{$\SU(6)$} & $T^2\subset\SU(3)\times\SU(3)\subset\SU(6)$\\
\hline
\multirow{2}*{$r=7$} & \multirow{2}*{$\varphi_4$} & \multirow{2}*{$\SU(8)$} & a finite group at least of order $128\subset\SU(4)\times\SU(4)$\\
\hline
\end{tabularx}\\
\end{center}

\quad II. $\fk=\begin{picture}(6,6)
\put(0,0){$B_r(r\geq2):$}
\put(50,0){$\circ$}
\put(53.5,2){\line(1,0){20}}
\put(73.5,0){$\circ$}
\put(78,2){\line(1,0){4}}
\put(84,1.5){$\ldots$}
%\put(50,1.5){$\ldots$}
\put(97,2){\line(1,0){6}}
\put(103,0){$\circ$}
\put(107,1){\line(1,0){20}}
\put(107,3){\line(1,0){20}}
%\put(65.3,1.3){\line(1,0){10}}
\put(126,0){$\bullet$}
\put(49,-6){$\alpha_1$}
\put(71.5,-6){$\alpha_2$}
\put(99,-6){$\alpha_{r-1}$}
\put(124,-6){$\alpha_r$}
\end{picture}
\quad\\$

\begin{center}
\begin{tabularx}{\textwidth}{|c|c|c|Y|}
\hline
rank $r$ &\hspace{5mm}The representation $\rho$ \hspace{5mm} &\hspace{10mm} The group $K$ \hspace{10mm} & $(H_{\rho})$\\
\hline
\multirow{3}*{$r\geq2$}& $\varphi_1$ & $\SO(2r+1)$ & $\SO(2r)$ \\
\cline{2-4}
& $S^2\varphi_1-\theta$ & $\SO(2r+1)$ & maximal $\mathbb{Z}_2$-tori\\
\cline{2-4}
& $\Lambda^2\varphi_1$ & $\SO(2r+1)$ & maximal tori \\
\hline
$r=2$ & $2\varphi_2$ & $\mathrm{Spin}(5)=\Sp(2)$ & $\Sp(1)$ \\
\hline
$r=3$ & $\varphi_3$ & $\mathrm{Spin}(7)$ & $G_2$ \\
\hline
\multirow{1}*{$r=4$} & \multirow{1}*{$\varphi_4$} & \multirow{1}*{$\mathrm{Spin}(9)$} & $\mathrm{Spin}(7)$ (with $\mathrm{Spin}(9)/\mathrm{Spin}(7)=S^{15}$)\\
\hline
\end{tabularx}\\
\end{center}
\quad\\

\quad III. $\fk=\begin{picture}(6,6)
\put(0,0){$C_r(r\geq3):$}
\put(50,0){$\bullet$}
\put(53.5,2){\line(1,0){20}}
\put(73.5,0){$\bullet$}
\put(78,2){\line(1,0){4}}
\put(84,1.5){$\ldots$}
%\put(50,1.5){$\ldots$}
\put(97,2){\line(1,0){6}}
\put(103,0){$\bullet$}
\put(107,1){\line(1,0){20}}
\put(107,3){\line(1,0){20}}
%\put(65.3,1.3){\line(1,0){10}}
\put(126,0){$\circ$}
\put(49,-6){$\alpha_1$}
\put(71.5,-6){$\alpha_2$}
\put(99,-6){$\alpha_{r-1}$}
\put(124,-6){$\alpha_r$}
\end{picture}
\quad\\$

\begin{center}
\begin{tabularx}{\textwidth}{|c|c|c|Y|}
\hline
rank $r$ &\hspace{5mm}The representation $\rho$ \hspace{5mm} &\hspace{10mm} The group $K$ \hspace{10mm} & $(H_{\rho})$\\
\hline
\multirow{3}*{$r\geq3$} & $2\varphi_1$ & $\Sp(r)$ & $\Sp(r-1)$ \\
\cline{2-4}
& $S^2\varphi_1$ & $\Sp(r)$ & maximal tori \\
\cline{2-4}
& \multirow{1}*{$\wedge^2\varphi_1-\theta$} & \multirow{1}*{$\Sp(r)/\mathbb{Z}_2$} & $(\Sp(1)\times\dots\times\Sp(1))/\mathbb{Z}_2$ ($r$ copies) \\
\hline
\multirow{1}*{$r=3$} & \multirow{1}*{$2\varphi_3$} & \multirow{1}*{$\Sp(3)$} & discrete subgroups isomorphic to $(\mathbb{Z}_2)^2$ \\
\hline
$r=4$ & $\varphi_4$ & $\Sp(4)$ & discrete subgroup of order $16$ \\
\hline
\end{tabularx}\\
\end{center}
\quad\\

\quad IV. $\fk=\begin{picture}(6,6)
\put(0,0){$D_r(r\geq4):$}
\put(50,0){$\circ$}
\put(53.5,2){\line(1,0){20}}
\put(73.5,0){$\circ$}
\put(78,2){\line(1,0){4}}
\put(84,1.5){$\ldots$}
%\put(50,1.5){$\ldots$}
\put(97,2){\line(1,0){6}}
\put(103,0){$\circ$}
\put(107,1){\line(2,-1){20}}
\put(107,3){\line(2,1){20}}
%\put(65.3,1.3){\line(1,0){10}}
\put(126,10){$\circ$}
\put(126,-11){$\circ$}
\put(49,-6){$\alpha_1$}
\put(71.5,-6){$\alpha_2$}
\put(97,-6){$\alpha_{r-2}$}
\put(124,-16){$\alpha_r$}
\put(124,6){$\alpha_{r-1}$}
\end{picture}
\quad\\$
\quad\\
%(There are automorphisms which interchange $\varphi_{r-1}, \varphi_r$ or $\varphi_1, \varphi_3, \varphi_4$ in case $r=4$)

\begin{tabularx}{\textwidth}{|c|c|c|Y|}
\hline
rank $r$ &\hspace{5mm}The representation $\rho$ \hspace{5mm} &\hspace{10mm} The group $K$ \hspace{10mm} & $(H_{\rho})$\\
\hline
\multirow{3}*{$r\geq4$}& $\varphi_1$ & $\SO(2r)$ & $\SO(2r-1)$ \\
\cline{2-4}
& $\varphi_2$ & $\SO(2r)$ & maximal tori\\
\cline{2-4}
& $S^2\varphi_1-\theta$ & $\SO(2r)$ & maximal $\mathbb{Z}_2$-tori \\
\hline
$r=5$ & $\varphi_4+\varphi_5$ & $\mathrm{Spin}(10)$ & $\SU(4)$ \\
\hline
$r=6$ & $2\varphi_5$ or $2\varphi_6$ & $\mathrm{Spin}(12)$ & $\SU(2)\times\SU(2)\times\SU(2)$ \\
\hline
$r=7$ & $\varphi_6+\varphi_7$ & $\mathrm{Spin}(14)$ & a finite group \\
\hline
$r=8$ & $\varphi_7$ or $\varphi_8$ & $\mathrm{Spin}(16)$ & a finite group \\
\hline
\end{tabularx}\\

\quad V. Exceptional Lie groups\\

\begin{tabularx}{\textwidth}{|Y|c|Y|}
\hline
The Lie algebra & The representation $\rho$ & $(H_{\rho})$ \\
\hline
\multirow{2}*{\begin{picture}(6,6)
\put(-70,0){$\G_2$:}
\put(-10,0){$\bullet$}
\put(-9,3){\line(1,0){30}}
\put(-9,1){\line(1,0){30}}
\put(20,0){$\circ$}
\put(-12,-5){$\alpha_1$}
\put(19,-5){$\alpha_2$}
\end{picture} }& $\varphi_1$ & maximal tori \\
\cline{2-3}
& $\varphi_2$ & $\SU(3)(\SU(3)/\G_2\approx S^6)$\\
\hline
\multirow{2}*{\begin{picture}(6,6)
\put(-70,0){$\F_4$:}
\put(-35,0){$\bullet$}
\put(-34,2){\line(1,0){30}}
\put(-5,0){$\bullet$}
\put(25,0){$\circ$}
\put(-4,3){\line(1,0){30}}
\put(-4,1){\line(1,0){30}}
\put(28,2){\line(1,0){28}}
\put(55,0){$\circ$}
\put(-37,-5){$\alpha_1$}
\put(-7,-5){$\alpha_2$}
\put(23,-5){$\alpha_3$}
\put(53,-5){$\alpha_4$}
\end{picture} }& $\varphi_1$ & maximal tori \\
\cline{2-3}
& $\varphi_4$(dim $\varphi_4=26$) & $\mathrm{Spin}(8)\subset\mathrm{Spin}(9)\subset \F_4$\\
\hline
\multirow{4}*{\begin{picture}(6,6)
\put(-70,5){$\E_6$:}
\put(-50,-10){$\circ$}
\put(-47,-8){\line(1,0){18}}
\put(-30,-10){$\circ$}
\put(-10,-10){$\circ$}
\put(-10,10){$\circ$}
\put(-7,-7){\line(0,1){18}}
\put(-27,-8){\line(1,0){18}}
\put(-7,-8){\line(1,0){18}}
\put(10,-10){$\circ$}
\put(14,-8){\line(1,0){17}}
%\put(34,-8){\line(1,0){17}}
\put(30,-10){$\circ$}
%\put(50,-10){$\circ$}
\put(-55,-15){$\alpha_5$}
\put(-32,-15){$\alpha_4$}
\put(-10,-15){$\alpha_3$}
\put(-2,7){$\alpha_6$}
\put(7,-15){$\alpha_2$}
\put(27,-15){$\alpha_1$}
\end{picture} }& \multirow{2}*{$\varphi_6$} & \multirow{2}*{maximal tori} \\
&&\\
\cline{2-3}
& \multirow{2}*{$\varphi_1+\varphi_5$} & \multirow{2}*{$\mathrm{Spin}(8)\subset \F_4\subset \E_6$}\\
&&\\
\cline{2-3}
\hline
\multirow{4}*{\begin{picture}(6,6)
\put(-70,5){$\E_7$:}
\put(-50,-10){$\circ$}
\put(-47,-8){\line(1,0){18}}
\put(-30,-10){$\circ$}
\put(-10,-10){$\circ$}
\put(-10,10){$\circ$}
\put(-7,-7){\line(0,1){18}}
\put(-27,-8){\line(1,0){18}}
\put(-7,-8){\line(1,0){18}}
\put(10,-10){$\circ$}
\put(14,-8){\line(1,0){17}}
\put(34,-8){\line(1,0){17}}
\put(30,-10){$\circ$}
\put(50,-10){$\circ$}
\put(-55,-15){$\alpha_6$}
\put(-32,-15){$\alpha_5$}
\put(-10,-15){$\alpha_4$}
\put(-2,7){$\alpha_7$}
\put(7,-15){$\alpha_3$}
\put(27,-15){$\alpha_2$}
\put(47,-15){$\alpha_1$}
\end{picture} }& \multirow{2}*{$\varphi_6$} & \multirow{2}*{maximal tori} \\
&&\\
\cline{2-3}
& \multirow{2}*{$2\varphi_1$} & \multirow{2}*{$\mathrm{Spin}(8)\subset \E_6\subset \E_7$}\\
&&\\
\cline{2-3}
\hline
\multirow{4}*{\begin{picture}(6,6)
\put(-70,5){$\E_8$:}
\put(-50,-10){$\circ$}
\put(-47,-8){\line(1,0){18}}
\put(-30,-10){$\circ$}
\put(-10,-10){$\circ$}
\put(-10,10){$\circ$}
\put(-7,-7){\line(0,1){18}}
\put(-27,-8){\line(1,0){18}}
\put(-7,-8){\line(1,0){18}}
\put(10,-10){$\circ$}
\put(14,-8){\line(1,0){17}}
\put(34,-8){\line(1,0){17}}
\put(54,-8){\line(1,0){17}}
\put(30,-10){$\circ$}
\put(50,-10){$\circ$}
\put(70,-10){$\circ$}
\put(-55,-15){$\alpha_7$}
\put(-32,-15){$\alpha_6$}
\put(-10,-15){$\alpha_5$}
\put(-2,11){$\alpha_8$}
\put(7,-15){$\alpha_4$}
\put(27,-15){$\alpha_3$}
\put(47,-15){$\alpha_2$}
\put(67,-15){$\alpha_1$}
\end{picture} }& \multirow{4}*{$\Ad(\E_8)$} & \multirow{4}*{maximal tori} \\
&&\\
&&\\
&&\\
\hline
\end{tabularx}\\
\quad\\
\end{footnotesize}

\end{document}